# GENERALIZATIONS OF A THEOREM ABOUT THE BINOMIAL COEFFICIENT

(*Dedicated to my friend Xuenong Cao*)

**SHAOHUA ZHANG**

School of Mathematics of Shandong University
P. R. China

The Key Lab of Cryptography Technology and Information Security
Ministry of Education
Shandong University
P. R. China
e-mail: shaohuazhang@mail.sdu.edu.cn

## Abstract

The object of this paper is to generalize a theorem on the binomial coefficient [4] to the case in an arithmetic progression. We will also give a slightly stronger result than Langevin's [2].

## 1 Introduction

As we know, $n!$ (the factorial of $n$) divides the product of $n$ consecutive positive integers. Namely, for any positive integer $m, n$, the fraction $\frac{(m+1)\times...\times(m+n)}{1\times 2\times...\times n}$ can be reduced to a positive integer. By [**4**], we know that $\frac{(m+1)\times...\times(m+n)}{1\times 2\times...\times n}$

------------------

This work was partially supported by the National Basic Research Program (973) of China (No. 2007CB807902) and the Natural Science Foundation of Shandong Province (No. Y2008G23).

Received April 20, 2009

may be written as $\dfrac{(m+1)\times...\times(m+n)}{1\times 2\times...\times n} = \prod_{i=1}^{n} a_i$, $a_i | m+i$, $a_i \in N$, $(a_i, a_j) = 1$,

$1 \leq i \neq j \leq n$. In this part, we will try to generalize this result to the case in an arithmetic progression. We will prove the following basic theorems:

**Theorem 1:** *For any positive integer $n$ and $a, b$ satisfying $(a,b) = 1$, $n!$ divides $\prod_{i=1}^{n}(a+bi)$ if and only if $(n!, b) = 1$.*

Based on the theorem 1, we generalize the theorem in [**4**] as following.

**Theorem 2:** *For any positive integer $n$ and $a, b$ satisfying $(a,b) = 1$, if $(n!, b) = 1$, then*

$\dfrac{(a+b)\times...\times(a+bn)}{1\times 2\times...\times n}$ *may be written as* $\dfrac{(a+b)\times...\times(a+bn)}{1\times 2\times...\times n} = \prod_{i=1}^{n} a_i$, $a_i | (a+bi)$,

$a_i \in N$, $(a_i, a_j) = 1$, $1 \leq i \neq j \leq n$.

Moreover, we generalize the theorem in [**5**] as follows.

**Theorem 3:** *For any positive integer $n$ and $a, b$ satisfying $(a,b) = 1$, if $(n!, b) = 1$ and*

$a + b > \prod_{p \leq n} p^{[\log_p^n]}$, *then* $\dfrac{(a+b)\times...\times(a+bn)}{1\times 2\times...\times n}$ *has the representation:*

$\dfrac{(a+b)\times...\times(a+bn)}{1\times 2\times...\times n} = \prod_{i=1}^{n} a_i$, $a_i | a+bi$, $a_i \in N$, $a_i > 1$, $(a_i, a_j) = 1$, $1 \leq i \neq j \leq n$.

*This implies that* $\prod_{i=1}^{n}(a+bi)$ *has $n$ distinct prime factors when* $(a,b) = 1$, $(n!, b) = 1$ *and*

$a + b > \prod_{p \leq n} p^{[\log_p^n]}$.

## 2 Proofs of Basic Theorems

For a positive integer $n > 1$, we denote by $p(n)$ the largest integer $e$ such that $p^e | n$, where $p$ is a prime factor of $n$. In this paper, we always denotes by $p$ a prime number.

**Proof of Theorem 1:** We write the factorial $n! = \prod_{i=1}^{r} p_i^{e_i}$, where $p_i$ satisfy $b \not\equiv 0 \pmod{p_i}$

($1 \leq i \leq r$). Thus, there must be positive integers $x_i, y_i$ such that

$$bx_i = p_i^{e_i} y_i + 1 \text{ for } 1 \leq i \leq r.$$

Hence, we have $(x_i a + x_i b)...(x_i a + x_i nb) \equiv (x_i a + 1)...(x_i a + n) \pmod{p_i^{e_i}}$. But, $n!$ divides

the product of $n$ consecutive positive integers. So $n! \mid (x_i a + 1)...(x_i a + n)$ and

$(x_i a + x_i b)...(x_i a + x_i nb) \equiv (x_i a + 1)...(x_i a + n) \equiv 0 \pmod{p_i^{e_i}}$. Note that $(x_i, p_i^{e_i}) = 1$.

Therefore, $(a+b)...(a+nb) \equiv 0 \pmod{p_i^{e_i}}$ and $n! \mid (a+b)...(a+nb)$. On the other hand,

if $n! \mid (a+b)...(a+nb)$, then, $(n!, b) = 1$ (otherwise $(a, b) > 1$). This completes the proof.

**Proof of Theorem 2:** Based on the methods in [5] or [4], in order to prove that Theorem 2 holds,

it is enough to prove that $p(\frac{(a+b) \times ... \times (a+bn)}{1 \times 2 \times ... \times n}) \leq t$ for any prime factor $p$ of $n!$, where $t = \max_{1 \leq i \leq n} \{p(a+bi)\}$. Let $p(n!) = e$. Thus, there must be positive integers $x, y$ such that $bx = p^e y + 1$ since $(n!, b) = 1$.

If $t < e$, then, we have $p((ax + xbi)) = p(ax + i + p^e yi) = p(ax + i)$. So

$$p((a+b)...(a+nb)) = p(\prod_{i=1}^{n}(ax + bxi)) = p(\prod_{i=1}^{n}(ax + i)).$$

By Lemma 1 in [4], we have $p\binom{ax+n}{n} \leq t$. Therefore, we have

$$p(\frac{(a+b) \times ... \times (a+bn)}{1 \times 2 \times ... \times n}) \leq t.$$

Now we consider the case $t \geq e$. If $e = 1$, then $n/2 < p \leq n$ since $p(n!) = e = 1$. So, at most two terms in the arithmetic progression $a+b, ..., a+nb$ can be divisible by $p$.

Note that these two terms can not be simultaneously divisible by $p^2$. Thus, $p(\frac{(a+b)\times...\times(a+bn)}{1\times 2\times...\times n}) \leq t$. If $e > 1$, then $n < p^e$ since $p(n!) = e$. So, at most one term in the arithmetic progression $a+b,...,a+nb$ can be divisible by $p^e$.

Let $p(a+br) = t = \max\limits_{1\leq i\leq n}\{p(a+bi)\}$. Note that $t \geq e$. So we have

$$p((ax+xbi)) = p(ax+i+p^e yi) = p(ax+i) \text{ for } i \neq r.$$

Thus, $p((a+b)...(a+nb)) = p(\prod_{i=1}^{n}(ax+bxi)) = p(\prod_{i\neq r}^{n}(ax+i)) + p(ax+bxr)$. Note that

$p(\prod_{i\neq r}^{n}(ax+i)) + p(ax+bxr) = p(\prod_{i=1}^{n}(ax+i)) + p(a+br) - p(ax+r)$. So, by Lemma 1 in

[**4**], $p(\frac{(a+b)\times...\times(a+bn)}{1\times 2\times...\times n}) = p(\prod_{i=1}^{n}(ax+i)) - p(n!) + p(a+br) - p(ax+r) \leq t$ holds

(since at most one term in the arithmetic progression $a+b,...,a+nb$ can be divisible by $p^e$, and, $p(ax+r) = \max\limits_{1\leq i\leq n}\{p(ax+i)\}$). Therefore, Theorem 2 holds.

**Proof of Theorem 3:** Since $a+b > \prod_{p\leq n} p^{[\log_p^n]}$, hence, there is a prime number $p_i$ such that

$p_i^{p_i(a+bi)} > n$ for any $1 \leq i \leq n$. Based on the methods in [**5**], in order to prove that Theorem 3

holds, it is enough to prove that $1 \leq p_i(\frac{(a+b)\times...\times(a+bn)}{1\times 2\times...\times n})$. Note that if $p_i^{p_i(a+bi)} > n$, then at most one term in the arithmetic progression $a+b,...,a+nb$ can be divisible by $p_i^{p_i(a+bi)}$.

Since $(n!, b) = 1$, hence there must be positive integers $x, y$ such that $bx = p_i^{p_i(a+bi)} y + 1$. So, we have $p_i(a+b) = p_i(ax+bx) = p_i(ax+1),..., p_i(a+(i-1)b) = p_i(ax+i-1)$;

$$p_i(a+ib) = p_i(ax+ibx) = p_i(ax+i+p_i^{p_i(a+bi)} yi);$$

$$p_i(a+(i+1)b) = p_i(ax+i+1),..., p_i(a+bn) = p_i(ax+n).$$

So, $p_i(\dfrac{(a+b) \times ... \times (a+bn)}{1 \times 2 \times ... \times n}) = \sum_{j=1}^{n} p_i(ax+j) + p_i(a+bi) - p_i(ax+i) - p_i(n!)$.

Note that $p_i(a+ib) = p_i(ax+i+p_i^{p_i(a+bi)}yi)$. So, $p_i(a+bi) \leq p_i(ax+i)$.

If $p_i(a+bi) = p_i(ax+i)$, by Lemma 2 in [**5**], then,

$$1 \leq p_i(\dfrac{(a+b) \times ... \times (a+bn)}{1 \times 2 \times ... \times n}).$$

If $p_i(a+bi) < p_i(ax+i)$, we write

$$p_i(\dfrac{(a+b) \times ... \times (a+bn)}{1 \times 2 \times ... \times n}) = \sum_{1 \leq j \neq i \leq n} p_i(ax+j) + p_i(a+bi) - p_i(n!).$$

Note that

$$\sum_{1 \leq j \neq i \leq n} p_i(ax+j) = \sum_{k=1}^{\infty}\left(\left[\dfrac{ax+i-1}{p_i^k}\right] - \left[\dfrac{ax}{p_i^k}\right] + \left[\dfrac{ax+n}{p_i^k}\right] - \left[\dfrac{ax+i}{p_i^k}\right]\right).$$

$$\sum_{k=1}^{\infty}\left(\left[\dfrac{ax+i-1}{p_i^k}\right] - \left[\dfrac{ax+i}{p_i^k}\right]\right)$$

$$= -p_i(a+bi) + \sum_{k=p_i(a+bi)+1}^{\infty}\left(\left[\dfrac{ax+i-1}{p_i^k}\right] - \left[\dfrac{ax+i}{p_i^k}\right]\right).$$

Based on the proof of Lemma 2 in [**5**], we have

$$\sum_{k=1}^{p_i(a+bi)}\left(-\left[\dfrac{ax}{p_i^k}\right] + \left[\dfrac{ax+n}{p_i^k}\right] - \left[\dfrac{n}{p_i^k}\right]\right) \geq 1.$$

Therefore, $p_i(\dfrac{(a+b) \times ... \times (a+bn)}{1 \times 2 \times ... \times n}) = \sum_{1 \leq j \neq i \leq n} p_i(ax+j) + p_i(a+bi) - p_i(n!)$

$$\geq 1 + \sum_{k=p_i(a+bi)+1}^{\infty}\left(\left[\dfrac{ax+i-1}{p_i^k}\right] - \left[\dfrac{ax}{p_i^k}\right] + \left[\dfrac{ax+n}{p_i^k}\right] - \left[\dfrac{ax+i}{p_i^k}\right] - \left[\dfrac{n}{p_i^k}\right]\right).$$

Note that

$$\left[\dfrac{n}{p_i^{p_i(a+bi)+1}}\right] = 0, \left[\dfrac{ax+i-1}{p_i^k}\right] \geq \left[\dfrac{ax}{p_i^k}\right], \left[\dfrac{ax+n}{p_i^k}\right] \geq \left[\dfrac{ax+i}{p_i^k}\right] \text{ since } 1 \leq i \leq n.$$

So, $\sum_{k=p_i(a+bi)+1}^{\infty}\left(\left[\dfrac{ax+i-1}{p_i^k}\right]-\left[\dfrac{ax}{p_i^k}\right]+\left[\dfrac{ax+n}{p_i^k}\right]-\left[\dfrac{ax+i}{p_i^k}\right]-\left[\dfrac{n}{p_i^k}\right]\right)\geq 0$ and this completes the proof of theorem 3.

## 3 Notes and Our Conclusions

In 1969, C. A. Grimm [1] conjectured that if $m+1,...,m+n$ are consecutive composite numbers, then there exist $n$ distinct prime numbers $p_1,...,p_n$ such that $m+j$ is divisible by $p_j$, $1\leq j\leq n$. Grimm's conjecture leads to the conjecture that there exists an injection $f$ from any finite set $\{a+b,...,a+bn\}$ of composite integers in arithmetic progression into the set of primes such that $f(a+bi)$ divides $a+bi$ for each $i=1,...,n$, where $a,b$ satisfy $(a,b)=1$. However, M. Langevin pointed out that Grimm's Conjecture cannot be extended to arithmetical progressions without a proviso: the numbers 12, 25, 38, 51, 64, 77, 90 belong to an arithmetic progression of ratio 13, but the number of distinct prime factors of $12\cdot 25\cdot 64\cdot 90$ is only 3. For details, see Michel Waldschmidt's *Open Diophantine Problems* [3]. Langevin [2] showed that the extension of Grimm's conjecture holds under the additional hypothesis that none of the integers $a+bi$, $i=1,...,n$, divides the least common multiple of the integers $1,...,n-1$. Note that the least common multiple of the integers $1,...,n-1$ is $\prod_{p\leq n-1} p^{[\log_p^{n-1}]}$. If none of the integers $a+bi$, $i=1,...,n$, divides $\prod_{p\leq n-1} p^{[\log_p^{n-1}]}$, then there must be a prime number $p_i$ such that $p_i^{p_i(a+bi)}>n-1$ for any $1\leq i\leq n$. Based on the proof of Theorem 3, it is easy to prove that for any positive integer $n$ and $a,b$ satisfying $(a,b)=1$, if $(n!,b)=1$ and none of the integers $a+bi$, $i=1,...,n$, divides $\prod_{p\leq n-1} p^{[\log_p^{n-1}]}$, then $\dfrac{(a+b)\times...\times(a+bn)}{1\times 2\times...\times n}$ has the representation: $\dfrac{(a+b)\times...\times(a+bn)}{1\times 2\times...\times n}=\prod_{i=1}^{n}a_i$, $a_i|a+bi$, $a_i\in N$, $a_i>1$, $(a_i,a_j)=1$,

$1 \leq i \neq j \leq n$. Moreover, by the following theorem 4, it is easy to prove that for every positive integer $n > 1$ and any positive integer $a, b$ satisfying $(a, b) = 1$, if none of the integers $a + bi$, $i = 1, ..., n$, divides the least common multiple of the integers $1, ..., n-1$, then $\frac{(a+b) \times ... \times (a+bn)}{1 \times 2 \times ... \times n}$ may be represented $\frac{(a+b) \times ... \times (a+bn)}{1 \times 2 \times ... \times n} = \frac{P}{Q}$, where $P = \prod_{i=1}^{n} a_i$, $a_i | a + bi$, $a_i \in N$, $a_i > 1$, $(a_i, a_j) = 1$, $1 \leq i \neq j \leq n$, and $Q$ satisfies $(Q, P) = 1$. Thus, we give a slightly stronger result than Langevin's [2]. Moreover, one could generalize basic theorems in Section 1 as follows:

**Theorem 4:** *For every positive integer $n > 1$ and any positive integer $a, b$ satisfying $(a, b) = 1$, $\frac{(a+b) \times ... \times (a+bn)}{1 \times 2 \times ... \times n}$ may be represented $\frac{(a+b) \times ... \times (a+bn)}{1 \times 2 \times ... \times n} = \frac{P}{Q}$, where $P = \prod_{i=1}^{n} a_i$, $a_i | a + bi$, $a_i \in N$, $(a_i, a_j) = 1$, $1 \leq i \neq j \leq n$, and $Q$ satisfies $(Q, P) = 1$.*

**Proof:** In order to prove Theorem 4, it is enough to prove that for any prime factor $p$ of $(a+b)...(a+nb)$, $p((a+b)...(a+nb)) \geq p(n!)$ holds. By the method of proof in Theorem 2, it is easy to prove that Theorem 4 holds. By Theorem 4, we know that for every positive integer $n > 1$ and any positive integer $a, b$ satisfying $(a, b) = 1$, if we write $\frac{(a+b) \times ... \times (a+bn)}{1 \times 2 \times ... \times n} = \frac{P}{Q}$ with $(Q, P) = 1$, then $P$ may be represented as follows $P = \prod_{i=1}^{n} a_i$, $a_i | a + bi$, $a_i \in N$, $(a_i, a_j) = 1$, $1 \leq i \neq j \leq n$. However, it is regret that the representation of $P$ is not unique. Although we have several methods for treating uniqueness, we do not know which treatment is the best. This problem will be given further consideration.

**Remark:** For an early version of this work, see [6]. However, the conjectures in [6] should be considered seriously, due to the fact that the theory is weak and the number evidence is still lacking. Although they are not disproved now, I would like to omit.

## Acknowledgements

I am thankful to the referee for his suggestions improving the presentation of the paper and also to my supervisor Professor Xiaoyun Wang for her valuable help and encouragement. Thank Doctor Jingguo Bi for his interesting discussion with me. Thank the Institute for Advanced Study in Tsinghua University for providing me with excellent conditions.

## References

[1]. C. A. Grimm, A conjecture on consecutive composite numbers, Amer. Math. Monthly, 76, (1969), 1126-1128.

[2]. M. Langevin, Plus grand facteur premier d'entiers en progression arithmétique. Séminaire Delange-Pisot-Poitou, 18e année: 1976/77, Théorie des nombres, Fasc. 1, Exp. No. 3, 7pp., Secrétariat Math., Paris, 1977.

[3]. Michel Waldschmidt, Open Diophantine Problems. Available at: http://arxiv.org/abs/math/0312440

[4]. Yongxing You and Shaohua Zhang, A new theorem on the binomial coefficient $\binom{m+n}{n}$, J. Math. (Wuhan, China) 23(2) (2003), 146-148.

[5]. Shaohua Zhang, A refinement of the function $g(m)$ on Grimm Conjecture, submitted to International Journal of Number Theory, submitted. Available at: http://arXiv.org/abs/0811.0966.

[6]. Shaohua Zhang, Generalizations of Several Theorems about the Binomial Coefficient and Grimm Conjecture. Available at: http://arxiv.org/abs/0901.3691v1